\documentclass[a4paper,reqno,twoside]{amsart}
\usepackage{amsmath,amsfonts,amssymb,enumerate,stmaryrd,verbatim,url}
\usepackage[english]{babel}

\theoremstyle{plain}
 \newtheorem{theorem}{Theorem}[section]
 \newtheorem{tthr}[theorem]{Theorem}
 \newtheorem{corr}[theorem]{Corollary}
 \newtheorem{lmr}[theorem]{Lemma}
 \newtheorem{prr}[theorem]{Proposition}
 \newtheorem{hyp}[theorem]{Conjecture}
\theoremstyle{definition}
 \newtheorem{defnr}[theorem]{Definition}
 \newtheorem{remnr}[theorem]{Remark}

\def\brm{\begin{remnr}}
\def\erm{\end{remnr}}
\def\bdr{\begin{defnr}}
\def\edr{\end{defnr}}
\def\bp{\begin{proof}}
\def\ep{\end{proof}}
\def\btr{\begin{tthr}}
\def\etr{\end{tthr}}
\def\bpr{\begin{prr}}
\def\epr{\end{prr}}
\def\bcr{\begin{corr}}
\def\ecr{\end{corr}}
\def\blr{\begin{lmr}}
\def\elr{\end{lmr}}
\def\beq{\begin{equation}}
\def\eeq{\end{equation}}
\def\bcs{\begin{cases}}
\def\ecs{\end{cases}}
\def\bhr{\begin{hyp}}
\def\ehr{\end{hyp}}

\def\en{\begin{enumerate}}
\def\ene{\end{enumerate}}

\def\d{\delta}
\def\dk{\delta_k}
\def\e{\varepsilon}
\def\ind{{\rm ind}}

\def\s{\prime}
\def\tl{\lhd}

\def\tlq{\trianglelefteqslant}

\def\w{\widehat}
\def\wg{\wedge}
\def\wdw{\wedge\dots\wedge}

\def\>{\geqslant}
\def\<{\leqslant}

\def\quat#1{\textquotedblleft #1\textquotedblright}

\title{COHOMOLOGY OF LIE ALGEBRAS OF POLYNOMIAL VECTOR FIELDS\\ ON THE LINE OVER FIELDS OF CHARACTERISTIC $2$}
\author{F.~V.~Weinstein}
\address
{Giacomettistrasse 33A, CH-3006 Bern, Switzerland.}
\email{felix.weinstein46@gmail.com}

\begin{document}

\begin{abstract}
For a field $\mathbb{F}$, let $L_k(\mathbb{F})$ be the Lie algebra
of derivations $f(t)\frac{d}{dt}$ of the polynomial ring
$\mathbb{F}[t]$, where $f(t)$ is a polynomial of degree $\>k$. For
any $k\geqslant -1$, we present a basis of the space of the cohomology with finite-dimensional support
of the Lie algebra $L_k(\mathbb{F})$ with coefficients in the
trivial module $\mathbb{F}$ for the case where
${\rm char}(\mathbb{F})=2$. The main result obtained is an analog of the famous
Goncharova's Theorem for the case ${\rm char}(\mathbb{F})=0$ and $k\>1$.
\end{abstract}

\maketitle

\markboth{F.~V.~WEINSTEIN}{COHOMOLOGY OF LIE ALGEBRAS OF POLYNOMIAL VECTOR FIELDS ON THE LINE}

\section*{\bf Introduction}

Let $W(\mathbb{F})$ be the vector space of polynomials in the
indeterminate $t$ over a field $\mathbb{F}$; we consider
${W(\mathbb{F})}$ as the Lie algebra with commutator
\[
\left[f_1(t),f_2(t)\right]=f_1(t)f^\prime_2(t)-f^\prime_1(t)f_2(t),
\]
where $f^\prime(t)$ denotes the derivative of the polynomial $f(t)$.
This algebra is called the \emph{Lie algebra of polynomial vector
fields on the line over $\mathbb{F}$}. The vectors $e_i=t^{i+1}$,
where $i\>-1$, form a basis of ${W(\mathbb{F})}$, and
\[
[e_a,e_b]=(b-a)e_{a+b}.
\]
The Lie algebra $W(\mathbb{F})$ contains a decreasing sequence of
the Lie subalgebras
\[
W(\mathbb{F})=L_{-1}(\mathbb{F})\supset
L_0(\mathbb{F})\supset L_1(\mathbb{F})\supset L_2(\mathbb{F})\supset\cdots,
\]
where $L_k(\mathbb{F})\subset W(\mathbb{F})$ is spanned by the vectors
$e_i$ with $i\>k$.

In this article we consider (co)homology of the Lie algebras $L_k(\mathbb{F})$ with
coefficients in the trivial module $\mathbb{F}$.
Since $L_k(\mathbb{F})$ is graded by integers $\>k$
(degrees), the corresponding chain Chevalley-Eilenberg complex
decomposes into a direct sum of the finite-dimensional subcomplexes
$C_*^{(n)}(L_k(\mathbb{F}))$,
enumerated by integers $n\>k$ (degrees of the chains).

Let us denote by $H^*_{(n)}(L_k(\mathbb{F}))$
the homology space of the dual complex $\mathrm{Hom}(C_*^{(n)}(L_k(\mathbb{F})),\mathbb{F})$
and define the cohomology of $L_k(\mathbb{F})$ by the formula
\[
H^*(L_k(\mathbb{F})):=\bigoplus_{n\>k}^\infty H^*_{(n)}(L_k(\mathbb{F})).
\]

Since for any $k\>-1$ the structure constants of the algebra
$L_k$ in the basis $e_i$ are integers, then to compute cohomology
of $L_k(\mathbb{F})$ over any field $\mathbb{F}$ it suffices to assume that $\mathbb{F}$
is a prime field.
That is, $\mathbb{F}=\mathbb{Q}$, the field of rational numbers, if ${\rm char}(\mathbb{F})=0$,
and $\mathbb{F}=\mathbb{Z}_p$, the field with $p$ elements, if ${\rm char}(\mathbb{F})=p>0$.

It turns out that the dimensions of the spaces
$H^q(L_k(\mathbb{Q}))$ are finite. It is easy to find them for
$k=-1,0$. For $k\>1$, these dimensions were predicted by D.~B.~Fuchs.
His formula was included in I.~M.~Gel$'$fand's talk \cite{MR0440631} as a conjecture.
It was proved by L.~V.~Goncharova in \cite{MR0339298}.

In article \cite{MR1254731} (see also \cite{W}) for any $k\>1$,
a special \quat{filtering} basis in the dual to the chain Chevalley-Eilenberg
complex of $L_k(\mathbb{Q})$ was offered. In this basis, the
coboundary operator acts simply enough and in particular allows to
obtain a transparent proof of Goncharova's result.

For $p>0$ and $q>0$ the spaces $H^q(L_k(\mathbb{Z}_p))$
in general are infinite--dimensional.
But the above definition implies that each of them is a direct sum of the
uniquely defined finite--dimensional subspaces:
\[
H^q(L_k(\mathbb{Z}_p))\simeq
\bigoplus_{n\>qk+\frac{q(q-1)}{2}}H^q_{(n)}(L_k(\mathbb{Z}_p)).
\]

The importance of Goncharova's result is explained in the book by D.~B.~Fuchs \cite{MR874337}.
In its Russian edition, as well as in private conversations, D.~B.~Fuchs
repeatedly asked \quat{what is the analog of Goncharova's theorem over fields of characteristic $p>0$ ?}.
Here I offer an answer for $p=2$.

In this article, for any $k\>1$, I build an analog of \quat{filtering} basis
for the Lie algebra $L_k(\mathbb{Z}_2)$
and apply it to the computation of space $H^*(L_k(\mathbb{Z}_2))$.
In \S\ref{Add}, I formulate a theorem
that explicitly describes a basis of $H^*(L_k(\mathbb{Z}_2))$ and compute the
numbers $\dim H^q_{(n)}(L_k(\mathbb{Z}_2))$. A detailed proof is
given only in the most interesting case $k=1$ in \S\ref{Coh1}.
For $k>1$, one can
prove the theorem from \S\ref{Add} by induction on $k$, similarly to what was done in
\cite{MR1254731} for $\mathbb{F}=\mathbb{Q}$.

The construction of a \quat{filtering} basis for $\mathbb{F}=\mathbb{Z}_2$ follows the
same scheme as for $\mathbb{F}=\mathbb{Q}$:
first, I build a special set of vectors, which spans the space $C^*(L_1(\mathbb{Z}_2))$. Their
linear independence follows from a combinatorial formula.
When $\mathbb{F}=\mathbb{Q}$ this formula is equivalent to the
Sylvester identity in the partitions theory (see
\cite{MR1634067} or \cite{W1}).
The corresponding formula in case of $\mathbb{F}=\mathbb{Z}_2$ leads to
another interesting identity for partitions,
which is established in \cite{W1} (see formula \eqref{Syl2} below).

This article is organized as follows. In \S\ref{sec_complex}, for
any $k\>1$, I introduce a complex, for which the cohomology will be
computed.

In \S\ref{p_Part} -- \S\ref{mult}, I consider only the case $k=1$.
The definitions related with integer partitions necessary to formulate the theorem on a
\quat{filtering} basis for $k=1$ are collected in \S\ref{p_Part}.

The \quat{filtering} basis theorem (Theorem \ref{T_Bas}) is formulated and proved in
\S\ref{sec_basis}.

In \S\ref{Coh1}, I use this result to build a basis of the bigraded
space $H^*(L_1(\mathbb{Z}_2))$ and compute the dimensions of its
homogeneous components (Theorem \ref{H1}).

In \S\ref{Add}, I formulate a theorem describing a basis of the
space $H^*(L_k(\mathbb{Z}_2))$ for any $k\>1$, thus generalizing the
result obtained for $k=1$ in \S\ref{Coh1}.

In \S\ref{mult}, I  compute, partly, the multiplication in
$H^*(L_1(\mathbb{Z}_2))$. In particular, I show that, as algebra, $H^*(L_1(\mathbb{Z}_2))$
is generated by $1$- and $2$-di\-men\-sional cohomological classes with a non-trivial product.
In addition I formulate a conjecture that completely describes the
multiplicative structure of $H^*(L_1(\mathbb{Z}_2))$.
Observe that the multiplication
in the algebra $H^*(L_1(\mathbb{Q}))$ is trivial (see \cite{W}, Remark 2.7).

In \S\ref{Coh0-1}, I use the results of \S\ref{Coh1}
to compute the cohomology of the Lie algebras $L_0(\mathbb{Z}_2)$ and $L_{-1}(\mathbb{Z}_2)$.
As an application, I offer an explicit description of the set
of equivalence classes of \emph{central extensions of
Lie algebra $W(\mathbb{Z}_2)=L_{-1}(\mathbb{Z}_2)$ with one-dimensional kernel},
i.e. the set of exact sequences of the Lie algebras
\[
0\longrightarrow\mathbb{Z}_2\longrightarrow\widetilde{W}(\mathbb{Z}_2)
\longrightarrow W(\mathbb{Z}_2)\longrightarrow 0
\]
considered up to the natural isomorphisms between such sequences (see \cite{Wei}, Sec.7.6).
\smallskip

\noindent \textbf{Notation:}
\smallskip

$|M|$ is the cardinality of the finite set $M$. \vspace{1mm}

$M_1\sqcup M_2$ is the disjoint union of the sets $M_1$ and $M_2$.\vspace{1mm}

For a finite set of integers $I$, set $I^-:=\min(I)$. \vspace{1mm}

For an integer $a$, let $\beta(a):=
\bcs 0&\text{if $a\equiv 0\mod 2$},\\
1&\text{if $a\equiv 1\mod 2$}. \ecs $ \vspace{1mm}

The symbol $L_k$ is used as a synonym of $L_k(\mathbb{Z}_2)$.\vspace{1mm}

In what follows all vector spaces are defined over the field $\mathbb{Z}_2$.

\section{\bf The complex $C^*(L_k)$}\label{sec_complex}

The \emph{homology of the Lie algebra $L_k$ with coefficients in the trivial
module $\mathbb{Z}_2$} is the homology of complex
\[
0\longleftarrow C_0(L_k)\overset{d}\longleftarrow C_1(L_k)
\overset{d}\longleftarrow\cdots\overset{d}\longleftarrow C_{q-1}(L_k)
\overset{d}\longleftarrow C_q(L_k)\longleftarrow\cdots,
\]
where $C_q(L_k):=\bigwedge^q L_k$
is the restricted (i.e., with finite support of each element) $q$th exterior power of the space $L_k$.
The vectors of the space $C_q(L_k)$ are called the \emph{$q$-chains of algebra $L_k$}.
Set ${C_*(L_k)=\bigoplus_{q=0}^\infty C_q(L_k)}$.

The set of $q$-chains
$e_{i_1}\wg e_{i_2}\wg\dots\wg e_{i_q}$,
where $i_1<i_2<\dots<i_q$, forms a basis of $C_q(L_k)$.
The action of the \emph{boundary operator} $d$
on such a chain is defined by the formula
\[
d(e_{i_1}\wg e_{i_2}\wdw e_{i_q})=\sum_{1\<a<b\<q}
(i_a+i_b)e_{i_a+i_b}\wg e_{i_1}\wdw\w e_{i_a}\wdw\w e_{i_b}\wdw e_{i_q}.
\]

This formula implies that the subspace
$C_*^{(n)}(L_k)\subset C_*(L_k)$, spanned by the chains ${e_{i_1}\wg e_{i_2}\wdw e_{i_q}}$
with $i_1+i_2+\dots+i_q=n$, is a finite-dimensional subcomplex,
and $C_*(L_k)=\bigoplus_{n\>0}C_*^{(n)}(L_k)$ is a direct sum of complexes.
(By definition, $C^{(n)}_0(L_k)=\{0\}$ for any appropriate integer $n$.)

Define the space of \emph{$q$-cochains} of algebra $L_k$ as
\[
C^q(L_k):=\bigoplus_{n\>0}C^q_{(n)}(L_k),\qquad\text{where}\qquad
C^q_{(n)}(L_k):=\text{\rm Hom}_\mathbb{F}\big(C_q^{(n)}(L_k),\mathbb{Z}_2\big).
\]
Set $C^*(L_k)=\bigoplus_{q=0}^\infty C^q(L_k)$.

On $C_*(L_k)$ let us introduce an Euclidean metric such that
the set of chains $e_{i_1}\wg e_{i_2}\wdw e_{i_q}$
with ${i_1<i_2<\dots<i_q}$ constitutes an orthonormnal basis.
For any $q$ and appropriate $n$ this metric defines
a unique isomorphism $C_q^{(n)}(L_k)\cong C^q_{(n)}(L_k)$ which allows us to
treat the $q$-chains of the algebra $L_k$ as its $q$-cochains.

Then the action of the \emph{coboundary operator} $\dk$,
i.e. the one dual to $d$, is expressed by the formulas:
\begin{gather}
\dk(e_i)=\sum\limits_{a+b=i;\ k\<a<b}\;(a+b)\,e_{a}\wg e_{b},\label{Bas_T_eq01}\\
\dk(e_{i_1}\wdw e_{i_q})=\sum_{1\<a\<q}\;
e_{i_1}\wdw\dk(e_{i_a})\wdw e_{i_q}.\label{Bas_T_eq02}
\end{gather}

\bdr For any integer $n\>0$, the cohomology space of the complex
\[
0\longrightarrow C_{(n)}^1(L_k)\overset{\d_k}\longrightarrow C_{(n)}^2(L_k)\overset{\d_k}\longrightarrow\cdots
\overset{\d_k}\longrightarrow C_{(n)}^{q-1}(L_k)\overset{\d_k}\longrightarrow
C_{(n)}^q(L_k)\overset{\d_k}\longrightarrow\cdots
\]
is denoted by $H^*_{(n)}(L_k)$.
\edr

\brm
Let $\mathcal{I}\subset L_{-1}$ be an ideal.
The adjoint action of $L_{-1}$ on $\mathcal{I}$
induces an action of $L_{-1}$ on ${C^*(\mathcal{I})}$,
which commutes with the coboundary operator.
Thus, $H^*(\mathcal{I})$ turns into an $L_{-1}$-module.
In our treatment of the cochains the action of $e_r\in L_{-1}$ on $C^*(\mathcal{I})$
is uniquely defined by the formula
\beq\label{act}
e_r(e_{i_1}\wdw e_{i_q})=\sum_{a=1}^q i_a\,e_{i_1}\wdw e_{i_a-r}\wdw e_{i_q}.
\eeq
\erm

\section{\bf Partitions}\label{p_Part}

This section gathers some material pertaining to partitions and is used throughout the article.

\bdr A \emph{partition} is a finite set of positive integers
$\langle i_1,i_2,\dots,i_q\rangle$, referred to the \emph{parts of the
partition} such that $i_1\<i_2\<\dots\<i_q$. A partition is said to
be \emph{strict} if $i_1<i_2<\dots<i_q$.

The numbers $\|I\|=i_1+i_2+\dots+i_q$ and $|I|=q$ are called the
\emph{degree} and \emph{length} of $I$, respectively.

A \emph{subpartition} of $I$ is an ordered subset $I^\prime\subset
I$. The \emph{union $I_1\sqcup I_2$ of partitions $I_1$ and $I_2$}
is the partition whose set of parts is the disjoint union of the
sets $I_1$ and $I_2$. \edr

\bdr
Let $I$ and $I^\s$ be distinct partitions such that $\|I\|=\|I^\s\|$. We write $I^\s\tl I$ if
either $|I^\s|<|I|$, or
\[
|I^\s|=|I|=q\quad\text{and}\quad i_1^\s+{\dots}+i_r^\s\<i_1+{\dots}+i_r\qquad\text{for any $r\in[1,q]$.}
\]
\edr

\bdr A \emph{marked partition} is a pair $\langle I;J\rangle$, where
$I$ is a partition and $J\subset I$. The elements of $J$ are called
the \emph{marked parts}. We identify each partition $I$ with the
marked partition $\langle I;\emptyset\rangle$. The numbers
\[
\left\|\langle I;J\rangle\right\|:=\|I\|,\qquad
\left|\langle I;J\rangle\right|:=|I|+|J|,\qquad\text{and}\qquad|I|
\]
are called the \emph{degree}, the \emph{length}, and the
\emph{reduced length} of $\langle I;J\rangle$, respectively.

A marked partition $\langle I;J\rangle$ is called \emph{strict} if $I$ is strict.

Define: $\langle I_1;J_1\rangle\sqcup\langle I_2;J_2\rangle:=
\big\langle I_1\sqcup I_2;J_1\sqcup J_2\big\rangle$.
\edr

Instead of explicitly indicating the set of marked parts,
we often underline these parts in $I$. For example,
$\langle 1,4,6,7\,;\,4,7\rangle=\langle 1,\underline{4},6,\underline{7}\rangle$.

\bdr
Let $\langle I;J\rangle$ and $\langle I^\s;J^\s\rangle$ be distinct marked partitions.
We write $\langle I^\s;J^\s\rangle\tl\langle I;J\rangle$ if
either $I^\prime\tl I$, or $I^\s=I$ and $J^\s\prec J$, where $\prec$
stands for the lexicographical order.
\edr

For example, $\langle\underline{5}\rangle\tl\langle 2,3\rangle,\,
\langle 3,\underline{9}\rangle\tl\langle 5,\underline{7}\rangle$ and
$\langle\underline{3},6\rangle\tl\langle 3,\underline{6}\rangle$.

One can readily show that the $\tl$ is a partial order on the set of
marked partitions. Below, we use the following

\blr[\cite{W}]\label{prmonBis}
If
$\langle I^\s;J^\s\rangle\tlq\langle I;J\rangle$ and
$\langle I^\s_1;J^\s_1\rangle\tl\langle I_1;J_1\rangle$,
then $\langle I^\s;J^\s\rangle\sqcup\langle I^\s_1;
J^\s_1\rangle\tl\langle I;J\rangle\sqcup\langle I_1;J_1\rangle$.
\elr

\bdr
A partition $I=\langle i_1,i_2,\dots i_q\rangle$ is called \emph{regular}
if $i_{m+1}-i_m\>2$ for any $m\in[1,q-1]$.

A \emph{dense} partition is a partition of the form $\langle a,a+2,a+4,\dots,a+2(q-1)\rangle$.

A \emph{special} partition is a dense partition of the form $\langle 1,3,\dots,2q-1\rangle$.

An \emph{even} or \emph{odd} partition is a partition all
parts of which are either even or odd, respectively. \edr

\bdr For any regular partition $I$, there is a unique decomposition
$I=I_1\sqcup I_2\sqcup\dots\sqcup I_s$, where ${I_1,I_2,\dots,I_s}$
are the dense subpartitions of the maximal possible length. This
decomposition is called the \emph{canonical decomposition of $I$}.
The dense partitions $I_1,I_2,\dots,I_s$ are said to be \emph{simple components} of $I$.

The minimal parts of the odd non-special simple components are called
the \emph{leading parts of $I$}. The quantity of leading parts is
denoted by $\ind(I)$. \edr

\bdr
A marked partition $\langle I;J\rangle$ is called \emph{regular}
if $I$ is regular and $J$ is a subset of the set of leading parts of $I$.
Otherwise it is called \emph{singular}.

For a regular marked partition $\langle I;J\rangle$, the
decomposition
\[
\langle I;J\rangle=\langle I_1;J_1\rangle\sqcup\langle I_2;
J_2\rangle\sqcup\dots\sqcup\langle I_s;J_s\rangle,
\]
where $I=I_1\sqcup I_2\sqcup\dots\sqcup I_s$ is the canonical decomposition,
is called the \emph{canonical decomposition of $\langle I;J\rangle$}.
The marked  $k$-partitions
$\langle I_1;J_1\rangle,\langle I_2;J_2\rangle,\dots,\langle I_s;J_s\rangle$
are called its \emph{simple components}.
\edr

\section{\bf A basis of $e$-monomials of the complex $C^*(L_1)$}\label{sec_basis}

\bdr\label{mark} Given $\langle I;J\rangle$, where $I=\langle
i_1,i_2,\dots i_q\rangle$, define
\[
e_{\langle I;J\rangle}=\widehat{e}_{i_1}\wedge\widehat{e}_{i_2}\wdw\widehat{e}_{i_q}\in
C^*(L_1),\quad\text{where}\qquad
\widehat{e}_{i_a}=
\bcs
e_{i_a}&\text{if $i_a\not\in J$},\\
\delta_1(e_{i_a})&\text{if $i_a\in J$}.
\ecs
\]
A cochain of the form $e_{\langle I;J\rangle}\in C^*(L_k)$ is called an \emph{$e$-monomial} if it is nonzero.
\edr

For example, $e_{\langle I;J\rangle}$ is an $e$-monomial for any regular marked
partition $\langle I;J\rangle$.

Keeping in mind the correspondence between marked partitions and
$e$-monomials, we apply the notions related to such partitions
(degree, length, regularity, order $\tl$, etc.), to $e$-monomials.
Lemma \ref{prmonBis} implies the following Lemma:

\blr\label{M_Part_S05} If $e_{\langle I;J\rangle}\tlq e_{\langle
I^\s;J^\s\rangle},\; e_{\langle I_1;J_1\rangle}\tl e_{\langle
I^\s_1;J^\s_1\rangle}$, and $e_{\langle I^\s;J^\s\rangle}\wg
e_{\langle I^\s_1;J^\s_1\rangle}\neq 0$, then
\[
e_{\langle I;J\rangle}\wg e_{\langle I_1;J_1\rangle}\tl e_{\langle
I^\s;J^\s\rangle}\wg e_{\langle I^\s_1;J^\s_1\rangle}.
\]
\elr

The next claim is the main result of this section:

\begin{theorem}\label{T_Bas} The set of regular $e$-monomials forms a basis of
the space $C^*(L_1)$. The decomposition of any singular $e$-monomial
$e_{\langle I;J\rangle}$ in this basis has the form
\beq\label{decompk}
e_{\langle I;J\rangle}=
\sum_{\langle I^\prime;J^\prime\rangle}e_{\langle I^\prime;J^\prime\rangle},
\eeq
where $\langle I^\prime;J^\prime\rangle\vartriangleleft\langle I;J\rangle$
for any $\langle I^\prime;J^\prime\rangle$
included in this sum.
\end{theorem}

\bp Let $D(n,q)$ be the set of partitions of degree $n$ and of length
$q$, and let $M(n,q)$ be the set of regular marked partitions of degree $n$
and of length $q$. Since $\dim C^q_{(n)}(L_1)=\left|D(n,q)\right|$, and assuming
the existence of decomposition \eqref{decompk}, linear independence
of the regular $e$-monomials follows from the combinatorial identity
\beq\label{Syl2}
\left|D(n,q)\right|=\left|M(n,q)\right|.
\eeq
This identity is proved in \cite{W1} (Theorem 2.1 for $\lambda=2$).

Therefore, it is sufficient to establish only the existence of the
decomposition \eqref{decompk}.

The set of singular $e$-monomials of reduced
length 2 is exhausted by the following $e$-monomials
\[
e_i\wg e_{i+1},\qquad
e_{2i}\wg\d_1(e_{2i\pm 1}),\qquad e_{2i-1}\wg\d_1(e_{2i\pm 1}),
\qquad\d_1(e_{2i-1})\wg\d_1(e_{2i+1}),
\]
where $i\>1$ is an appropriate integer. For this set of
$e$-monomials of degree $n$, the decomposition \eqref{decompk}
directly follows from the easily checked identities
\begin{gather}
\sum_{a+b=n,\,1\<a<b}e_a\wg e_b=\d_1(e_n)\qquad\text{if $n\equiv 1\mod 2$},\label{dim21}\\
\sum_{a+b=n,\,a,b\>1}e_a\wg\d_1(e_b)=0,\label{dim22}\\[1mm]
\sum_{a+b=n,\,1\<a<b}\d_1(e_a)\wg\d_1(e_b)=0
\qquad\text{if $n\not\equiv 2\mod 4$}\label{dim23}.
\end{gather}
Using these identities we can express any singular submonomial of reduced length $\<2$ in
$e_{\langle I;J\rangle}$ as a linear combination of the regular
ones. As a result, we present the monomial $e_{\langle I;J\rangle}$
as a linear combination of $e$-monomials
$e_{\langle I^\prime;J^\prime\rangle}$ such that
$e_{\langle I^\prime;J^\prime\rangle}\tl e_{\langle I;J\rangle}$ as it follows
from Lemma \ref{M_Part_S05}.

Let us now apply the same procedure to each singular $e$-monomial
$e_{\langle I^\prime;J^\prime\rangle}$ of the obtained linear
combination. Since the set of $e$-monomials of fixed degree is
finite, a finite number of such iterations leads to the
decomposition \eqref{decompk}. \ep

\section{\bf Computing the space $H^*(L_1)$}\label{Coh1}

Let $M_0$ be the set of marked regular partitions. For $r\>1$,
define
\[
M_r:=M_{r-1}\smallsetminus\left\{\text{the set of maximal partitions (with respect to order $\tl$)
in $M_{r-1}$}\right\}.
\]
In view of Theorem \ref{T_Bas}, the sequence of partially ordered
sets $M_0\supset M_1\supset M_2\supset\cdots$ induces a filtration
of the vector spaces \beq\label{filt} C^*(L_1)=C^*( M_0)\supset C^*(
M_1)\supset C^*(M_2)\supset\cdots, \eeq where $C^*( M_r)\subset
C^*(L_1)$ is the subspace spanned by regular
$e$-monomials $e_{\langle I;J\rangle}$ such that
${\langle I;J\rangle\in M_r}$.

For the canonical decomposition $\langle I;J\rangle=\langle
I_1;J_1\rangle\sqcup\dots\sqcup\langle I_s;J_s\rangle$, formula
\eqref{Bas_T_eq02} implies that
\beq\label{Bas_T_eq07}
\d_1\left(e_{\langle I;J\rangle}\right)=\sum_{1\<a\<s} e_{\langle
I_1;J_1\rangle}\wdw\d_1\left(e_{\langle I_a;J_a\rangle}\right)\wdw
e_{\langle I_s;J_s\rangle}.
\eeq

This expression, Theorem \ref{T_Bas}, and Lemma \ref{M_Part_S05}
show that \eqref{filt} is a filtration of the complex $C^*(L_1)$.

\bdr
For any $x,y\in C^*( M_r)$, we write $x\approx y$ if $x-y\in C^*( M_{r+1})$.
\edr

\blr\label{delta1}
Let $I$ be a dense partition. Then
\[
\d_1(e_I)\approx
\bcs
0&\text{if $I$ is special or even},\\
\beta(|I|)\;e_{\langle I;I^-\rangle}&\text{otherwise}.
\ecs
\]
\elr

\bp
For any special or even partition $I$, it is clear since $\d_1(e_I)=0$ for such a partition.

Formula \eqref{dim22} implies that
$e_{a+2r}\wg\d_1(e_{a+2(r+1)})\approx\d_1(e_{a+2r})\wg e_{a+2(r+1)}$
for any $r\>0$. Therefore, for the remaining partitions,
Lemma \ref{delta1} follows from formula \eqref{Bas_T_eq02} and Lemma \ref{M_Part_S05}.
\ep

\blr\label{delta2}
Let $I$ be a dense odd non-special partition of even length.
Then there exist cocycles $\e_I$ and $\e_{\langle I;I^-\rangle}$ of the complex $C^*(L_1)$ such that
\[
\e_I\approx e_I\qquad\text{and}\qquad \e_{\langle I;I^-\rangle}\approx e_{\langle I;I^-\rangle}.
\]
\elr

\bp
Formulas \eqref{dim22} and \eqref{dim23} imply that, for $a\>3$  odd, the cochains
\beq\label{co2}
\e_{\langle a,a+2\rangle}:=\sum_{r=0}^{\frac{a-1}{2}}e_{a-2r}\wg e_{a+2r+2},\qquad
\e_{\langle\underline{a},a+2\rangle}:=\sum_{r=0}^{\frac{a-1}{2}}\d_1(e_{a-2r})\wg e_{a+2r+2}
\eeq
are cocycles of the complex $C^*(L_1)$.
Therefore, for $I=\langle a,a+2,\dots a+2(q-1)\rangle$, $q$ even, and $a\>3$ odd,
formula \eqref{Bas_T_eq07} implies that the following cochains are cocycles as well:
\begin{align*}
\e_I:=\;&\e_{\langle a,a+2\rangle}\wg
\e_{\langle a+4,a+6\rangle}\wg\dots\wg\e_{\langle a+2(q-2),a+2(q-1)\rangle},\\
\e_{\langle I;I^-\rangle}:=\;&
\e_{\langle\underline{a},a+2\rangle}\wg\e_{\langle a+4,a+6\rangle}
\wg\dots\wg\e_{\langle a+2(q-2),a+2(q-1)\rangle}.
\end{align*}
The formulas $\e_I\approx e_I$ and $\e_{\langle I;I^-\rangle}\approx e_{\langle I;I^-\rangle}$
follow from Theorem \ref{T_Bas} and Lemma \ref{M_Part_S05}.
\ep
For dense partitions $I$ not mentioned in Lemma \ref{delta2}, define
\[
\e_I:=e_I,\qquad\text{and}\qquad
\e_{\langle I;I^-\rangle}:=\d_1(e_I)\quad
\text{if $I$ is odd, non-special, and $\beta(|I|)=1$}.
\]
These formulas, together with the formulas from Lemma \ref{delta2}, define cochains
$\e_I$ and $\e_{\langle I;I^-\rangle}$ for any dense regular partition $I$.
Cochains of such a form are said to be the \emph{simple $\e$-monomials}.

\bdr
Let
$\langle I;J\rangle=\langle I_1;J_1\rangle\sqcup\dots\sqcup\langle I_s;J_s\rangle$
be the canonical decomposition of a regular partition.
Define
$\e_{\langle I;J\rangle}:=\e_{\langle I_1;J_1\rangle}\wg\dots\wg\e_{\langle I_s;J_s\rangle}$.
The cochain $\e_{\langle I;J\rangle}$
is called the \emph{$\e$-monomial} corresponding to ${\langle I;J\rangle}$.
Its simple $\wg$-factors are called the \emph{simple components} of the
$\e$-monomial $\e_{\langle I;J\rangle}$.
\edr

Since the matrix of passage from the set of $\e$-monomials to the basis of regular $e$-monomials
is a square lower triangle matrix with units on the main diagonal,
Theorem \ref{T_Bas} and Lemmas \ref{M_Part_S05}, \ref{delta1}, and \ref{delta2}
imply the following result:

\btr\label{eposBas}
The set of $\e$-monomials is a basis of the complex $C^*(L_1)$.
For a simple $\e$-monomial ${\e_{\langle I;J\rangle}}$, we have:
\[
\d_1\left(\e_{\langle I;J\rangle}\right)=
\bcs
\e_{\langle I;I^-\rangle}&\text{if $I$ is an odd non-special partition,
$J=\emptyset$, and $\beta(|I|)=1$},\\
0&\text{otherwise}.
\ecs
\]
\etr

From this theorem, formula \eqref{Bas_T_eq07}, and the definition of
the tensor product of complexes, we obtain

\bcr\label{deps}
For a regular partition $I$, let $T^*(I)\subset C^*(L_1)$ be the linear span
of the set of $\e$-monomials $\e_{\langle I;J\rangle}$.
Then $T^*(I)$ is a subcomplex in $C^*(L_1)$ and
\[
C^*(L_1)=\bigoplus_{I}\;T^*(I),
\]
where $I$ runs over the set of regular partitions.
For the canonical decomposition
$I=I_1\sqcup I_2\sqcup\dots\sqcup I_s$, we have an isomorphism of complexes
\[
T^*(I)\simeq T^*(I_1)\otimes T^*(I_2)\otimes\cdots\otimes T^*(I_s).
\]
\ecr

\btr\label{H1}
Let $R(n)$ be the set of the regular partitions of
degree $n$, for which each simple component is either special, or
has an even degree.

Then any $\e$-monomial $\e_{\langle I;J\rangle}$,
where $I\in R(n)$, is a nonzero cocycle of the complex $C^*(L_1)$.
Cohomological classes of these cocycles form a basis of the space $H^*_{(n)}(L_1)$.
Moreover, the classes with ${|I|+|J|=q}$
form a basis of the space $H^q_{(n)}(L_1)$.
In particular,
\[
\sum_{q=1}^\infty\,\dim H^q_{(n)}(L_1)\,t^q=\sum_{I\in R(n)}\;(1+t)^{\ind(I)}\;t^{|I|}.
\]
\etr

\bp
Theorem \ref{eposBas} implies that, for a dense partition $I$, we have
\[
H^*(T^*(I))=
\bcs
0&\text{if $I$ is odd non-special partition and $\beta(|I|)=1$},\\
T^*(I)&\text{otherwise}.
\ecs
\]
Therefore Theorem \ref{H1} follows from
Corollary \ref{deps} and the K\"unneth isomorphism
\[
H^*(T^*(I))\simeq H^*(T^*(I_1))\otimes H^*(T^*(I_2))\otimes\cdots\otimes H^*(T^*(I_s)).
\qedhere
\]
\ep

For example,
$R(12)=\{\langle 12\rangle,\langle 2,10\rangle,\langle 4,8\rangle,
\langle 5,7\rangle,\langle 1,3,8\rangle,\langle 2,4,6\rangle\}$.
If $I\in R(12)$ and $I\neq\langle 5,7\rangle$, then $\ind(I)=0$.
Since $\ind(\langle 5,7\rangle)=1$, we have
\[
\sum_{q=1}^\infty\dim H^q_{(12)}(L_1)\,t^q=t+3t^2+3t^3.
\]

\section{\bf A basis of the space $H^*(L_k)$ for $k\>1$}\label{Add}

\bdr\label{defk}
A \emph{$k$-partition} is a pair $(k;I)$, where $I$ is a
partition and $I^-\>k$.
\edr
We write $k$-partitions as usual
partitions, emphasizing that we only consider $k$-partitions. For
example, one may consider $\langle 2,7\rangle$ as either a
1-partition or a 2-partition; these objects are different.

\bdr
A regular $k$-partition $I=\langle i_1,i_2,\dots i_q\rangle$ is called \emph{special}
if $i_q<2(k+q-1)$.
\edr

It is easy to check that the quantity of the special $k$-partitions
of length $q$ is equal to $\binom{q+k-1}{k-1}$.

\bdr
A regular $k$-partition is called \emph{simple} if it is either special, or dense.
\edr
\bdr
For any regular $k$-partition $I$, there exists a unique decomposition
${I=I_1\sqcup I_2\sqcup\dots\sqcup I_s}$,
where $I_1,I_2,\dots,I_s$ are simple $k$-partitions of the maximal possible length.
This decomposition is called the \emph{canonical decomposition of $I$}.
The partitions $I_1,I_2,\dots,I_s$ are called the \emph{simple components of $I$}.

The minimal parts of the odd non-special simple components of $I$
are called the \emph{leading parts of $I$}.
The quantity of them is denoted by $\ind_k(I)$.
\edr

For example, for $I=\langle 3,5,9,13,15,18\rangle$
and $k=1,2,3$, the canonical decompositions are
\[
I=
\bcs
\langle 3,5\rangle\sqcup\langle 9\rangle\sqcup\langle 13,15\rangle
\sqcup\langle 18\rangle&\text{if $k=1$},\\
\langle 3,5\rangle\sqcup\langle 9\rangle\sqcup\langle 13,15\rangle
\sqcup\langle 18\rangle&\text{if $k=2$},\\
\langle 3,5,9\rangle\sqcup\langle 13,15\rangle\sqcup
\langle 18\rangle&\text{if $k=3$}.
\ecs
\]
Observe that decompositions for $k=1$ and $k=2$
are distinct as k-partitions (see Def.\ref{defk}).
The definition of $\ind_k(I)$ implies that
$\ind_1(I)=3$, $\ind_2(I)=2$, and $\ind_3(I)=1$, respectively.

\bdr
We say that $\langle I;J\rangle$ is a
\emph{$k$-partition} if $I$ is a $k$-partition.
A marked $k$-partition $\langle I;J\rangle$ is called \emph{regular}
if $I$ is regular and $J$ is a subset of the set of the leading parts of $I$.
\edr

\btr\label{Hk0}
Let $R_k(n)$ be the set of $k$-regular partitions of degree $n$,
each simple component of which is either special, or of even degree.

For any regular marked $k$-partition $\langle I;J\rangle$, where $I\in R_k(n)$,
one can uniquely define a nonzero cocycle $\e_{\langle I;J\rangle}\in C^*_{(n)}(L_k)$.
Cohomological classes of these cocycles form a basis of the space $H^*_{(n)}(L_k)$.
Moreover, the classes  with $|I|+|J|=q$ form a basis of the space $H^q_{(n)}(L_k)$.
In particular,
\[
\sum_{q=1}^\infty\,\dim H^q_{(n)}(L_k)\,t^q=\sum_{I\in R_k(n)}\;(1+t)^{\ind_k(I)}\;t^{|I|}.
\]
\etr

\section{\bf On the multiplicative structure of the algebra  $H^*(L_1)$}\label{mult}

The exterior product of cochains in the
complex $C^*(L_1)$ induces a multiplication that turns ${H^*(L_1)}$ into an
algebra. Theorem \ref{H1} and formulas \eqref{co2} imply that
$\e$-monomials $e_1$ and
\begin{gather*}
x(i):=e_{2i},\qquad
y(i):=\sum_{r=0}^{i-1}e_{2i-2r-1}\wg e_{2i+2r+1},\qquad\text{where $i\>1$},\\
z(i):=\sum_{r=0}^{i-2}\d_1(e_{2i-2r-1})\wg e_{2i+2r+1},\qquad\text{where $i\>2$},
\end{gather*}
are cocycles which represent nonzero cohomological classes of $L_1$.
Let us denote these classes by ${e,\,x_i,\,y_i}$, and $z_i$, respectively.
Theorem \ref{H1} implies that they multiplicatively generate the algebra $H^*(L_1)$.

\blr\label{pi1}
In the algebra $H^*(L_1)$, we have $e^2=0,\,x^2_i=y^2_i=0$ for all $i\>1$ and
\begin{gather*}
e\cdot x_1=e\cdot y_1=0,\label{f1}\\
z_i=\sum_{a=1}^{i-1}\,x_{2a}\cdot y_{i-a}\qquad\text{for\, $i\>2$},\label{f2}\\
\sum_{a=0}^{i-1}\,x_{2a+1}\cdot y_{i-a}=0,\qquad
\sum_{a=0}^{i-1}\,y_{i-a}\cdot y_{i+a+1}=0 \qquad\text{for\, $i\>1$}\label{f3}.
\end{gather*}
In particular, the algebra $H^*(L_1)$ is multiplicatively
generated by the classes $e,\,x_i$, and $y_i$ for all ${i\>1}$.
\elr

\bp
Clearly, $e\wg x_1=e\wg y_1=0$. The remaining formulas follow
from the directly checked relations
\begin{gather*}
z(i)=\sum_{a=1}^{i-1}x(2a)\wg y(i-a)+
\d_1\left(\,\sum_{m=0}^{\left\lfloor\frac{i-2}{2}\right\rfloor}
e_{2i-4m-3}\wg e_{2i+4m+3}\right)\qquad\text{for\, $i\>2$},\\[1mm]
\sum_{a=0}^{i-1}x(2a+1)\wg y(i-a)=
\d_1\left(\,\sum\limits_{m=1}^{\left\lfloor\frac{i+1}{2}\right\rfloor}
e_{4m-1-2\beta(i)}\wg e_{4(i-m)+3+2\beta(i)}\right)\qquad\text{for\, $i\>1$},\\
\sum_{a=0}^{i-1}y(i-a)\wg y(i+a+1)=0\qquad\text{for\, $i\>1$}.
\end{gather*}
In addition, Lemma \ref{M_Part_S05} implies that
\[
e_1\wg e_3\wg\dots\wg e_{2q-1}\approx
\bcs\label{cor2}
y(1)\wg y(3)\wg\dots\wg y(q-1)&\text{if $\beta(q)=0$},\\
e_1\wg y(2)\wg y(4)\wg\dots\wg y(q-1)&\text{if $\beta(q)=1$}.
\ecs
\]
To finish the proof, it remains to apply Theorems \ref{H1} and \ref{T_Bas}.
\ep

Let $\mathbb{P}[E,X,Y]$ be the exterior algebra generated over
$\mathbb{Z}_2$ by $E,\,X_i$, and $Y_i$ for all $i\>1$.
Consider $\mathbb{P}[E,X,Y]$ as a bigraded algebra having defined the
bigrading as follows:
\[
\deg(E)=(1,1),\qquad\deg(X_i)=(1,2i),\qquad\deg(Y_i)=(2,4i).
\]

\bhr\label{hypm}
Let $A$ be a homogeneous ideal in $\mathbb{P}[E,X,Y]$ with generators
\[
E\wg X_1,\qquad E\wg Y_1,\qquad
\sum_{a=0}^{i-1}\,X_{2a+1}\wg Y_{i-a},\qquad
\sum_{a=0}^{i-1}\,Y_{i-a}\wg Y_{i+a+1},\qquad\text{where\, $i\>1$}.
\]
Then we have the exact sequence of bigraded algebras
\[
0\longrightarrow A\longrightarrow\mathbb{P}[E,X,Y]
\overset{\pi}\longrightarrow H^*(L_1)\longrightarrow 0,
\]
where $\pi(E)=e,\,\pi(X_i)=x_i,\,\pi(Y_i)=y_i$ for all $i\>1$.
\ehr

\brm
Thanks to Theorem \ref{H1} one can reduce this conjecture
to a combinatorial question about integer partitions. Namely,
let $P(n,q)$ be the set of pairs $(K,L)$ of partitions such that
\en
\item
$K=\langle k_1,k_2,\dots k_a\rangle$ is a strict partition and
$L=\langle l_1,l_2,\dots l_b\rangle$ is a regular partition.
\item
$a+2b=q$ and $2k_1+2k_2+\dots+2k_a+4l_1+4l_2+\dots+4l_b=n$.
\item
$|k_i-l_j|\>2$ for all $i\in[1,a]$ such that $\beta(k_i)=1$, and for all $j\in[1,b]$.
\ene
On the other hand, let $M_0(n,q)$ be
the set of marked regular partitions $\langle I;J\rangle$ of
degree $n$, of length $q$, and such that the degree of any simple component of $I$ is even.
Then it is easy to see that Conjecture \ref{hypm} follows from the conjectural
identity $|P(n,q)|=|M_0(n,q)|$.
\erm

\section{\bf Computing the spaces $H^*(L_0)$ and $H^*(L_{-1})$}\label{Coh0-1}

In this section, $k=0$ or $k=-1$.
Since the cohomology of $L_1$ are now known and in both the cases $L_1\subset L_k$ is
an ideal, to compute $H^*(L_k)$
one could use the corresponding Hochischid-Serre spectral sequence.
But we prefer a direct computation, which agrees with
the computations from \S\ref{Coh1}.

Let $A^q_{(n)}(k)$ be the vector subspace of $C^q_{(n)}(L_k)$ spanned by
the cochains $c=e_{i_1}\wg e_{i_2}\wg\dots\wg e_{i_q}$ such that
$e_0\wg c\neq 0$. Obviously,
\[
A^q_{(n)}(0)=C^q_{(n)}(L_1),\qquad
A^q_{(n)}(-1)=
\bcs
\mathbb{Z}_2\,e_n&\text{if $q=1$ and $n\neq 0$},\\
e_{-1}\wg C^{q-1}_{(n+1)}(L_1)\bigoplus C^q_{(n)}(L_1)&\text{if $q>1$ and $n\>0$},\\
\ecs
\]
Then
\[
C^q_{(n)}(L_k)=
\bcs
\mathbb{Z}_2\,e_n&\text{if $q=1$},\\
e_0\wg A^{q-1}_{(n)}(k)\bigoplus A^q_{(n)}(k)&\text{if $q>1$}.
\ecs
\]

Since $e_0$ is a cocycle of the complex $C^*(L_k)$,
the space $e_0\wg A^*_{(n)}(k)$ is a subcomplex of ${C^*_{(n)}(L_k)}$.
Using the natural projection, let us identify the space $A^*_{(n)}(k)$ with the space of
the factor-complex $C^*_{(n)}(L_k)/\left(e_0\wg A^*_{(n)}(k)\right)$
and transfer its differential, denoted by $\delta$, to $A^*_{(n)}(k)$.
Then
\[
H^q(e_0\wg A^*_{(n)}(k))\cong e_0\wg H^{q-1}( A^*_{(n)}(k)).
\]
The exact sequence of complexes
\[
0\longrightarrow e_0\wg A^*_{(n)}(k)\longrightarrow C^*_{(n)}(L_k)\longrightarrow A^*_{(n)}(k)\longrightarrow 0,
\]
induces the exact cohomology sequence
\[
{\cdots}\longrightarrow H^{q-1}_{(n)}(A^*(k))\overset{b}\longrightarrow e_0\wg H^{q-1}_{(n)}(A^*(k))\longrightarrow H^q_{(n)}(L_k)
\longrightarrow H^q_{(n)}(A^*(k))\overset{b}\longrightarrow e_0\wg H^q_{(n)}(A^*(k))\longrightarrow\cdots,
\]
where $b$ is the Bockstein homomorphism.
Its definition implies that, for any $h\in H^*_{(n)}(A^*(k))$, we have
\[
b(h)=ne_0\wg h
\]
Therefore, for integer $l\>0$ if $n=2l+1$ then $H^q_{(n)}(L_k)=\{0\}$, whereas
\beq\label{CH}
H^q_{(2l)}(L_k)\cong
\bcs
\mathbb{Z}_2\,e_{2l}&\text{if $q=1$},\\
e_0\wg H^{q-1}\big(A^*_{(2l)}(k)\big)\bigoplus H^q\big(A^*_{(2l)}(k)\big)&\text{if $q>1$}.
\ecs
\eeq

For $k=0$, this implies the following result:
\btr
We have, $H^q_{(2l+1)}(L_0)=\{0\}$, whereas
\[
H^q_{(2l)}(L_0)\cong
\bcs
\mathbb{Z}_2\,e_{2l}&\text{if $q=1$},\\
e_0\wg H^{q-1}_{(2l)}(L_1)\bigoplus H^q_{(2l)}(L_1)&\text{if $q>1$}.
\ecs
\]
\etr

Let now $k=-1$. Since $\delta(e_{-1})=0$, we see that
$e_{-1}\wg C^*_{(n+1)}(L_1)\subset A^*_{(n)}(-1)$
is a subcomplex and $H^q(e_{-1}\wg C^*_{(n+1)}(L_1))\cong e_{-1}\wg H^q_{(n+1)}(L_1)$.
The exact sequence of complexes
\[
0\longrightarrow e_{-1}\wg C^*_{(n+1)}(L_1)
\longrightarrow A^*_{(n)}(-1)\longrightarrow C^*_{(n)}(L_1)\longrightarrow 0,
\]
induces the exact cohomology sequence
\[
{\cdots}\longrightarrow H^{q-1}_{(n)}(L_1)
\overset{b}\longrightarrow e_{-1}\wg H^{q-1}_{(n+1)}(L_1)\longrightarrow H^q(A^*_{(n)}(-1))
\longrightarrow H^q_{(n)}(L_1)\overset{b}\longrightarrow e_{-1}\wg H^q_{(n+1)}(L_1)\longrightarrow\cdots.
\]
The definition of $b$ implies that for any $h\in H^*_{(n)}(L_1)$ we have
\[
b(h)=e_{-1}\wg e_{-1}(h),
\]
where $e_{-1}(h)$ denotes the action of $e_{-1}$ on $h\in H^*_{(n)}(L_1)$ (see formula \eqref{act}).
It is subject to a direct verification that for $a\>1$ odd, we have
\[
e_{-1}(\e_{a,a+2})=\d_1(e_{2a+3})\qquad\text{and}\qquad
e_{-1}(\e_{\underline{a},a+2})=\d_1(e_a\wg e_{a+3}).
\]
Since in addition $e_{-1}(e_i)=0$ for any $i\>2$ even,
Theorem \ref{H1} implies that $b(h)=0$. Thus,
\[
H^q(A^*_{(n)}(-1))\cong e_{-1}\wg H^{q-1}_{(n+1)}(L_1)\bigoplus H^q_{(n)}(L_1).
\]
Therefore, the isomorphism \eqref{CH} implies the following result:
\btr
We have $H^q_{(2l+1)}(L_{-1})=\{0\}$, whereas
\[
H^q_{(2l)}(L_{-1})\cong
\bcs
\mathbb{Z}_2\,e_{2l}&\text{if $q=1$},\\
e_{-1}\wg e_0\wg H^{q-2}_{(2l+1)}(L_1)\bigoplus e_0\wg H^{q-1}_{(2l)}(L_1)
\bigoplus e_{-1}\wg H^{q-1}_{(2l+1)}(L_1)\bigoplus H^q_{(2l)}(L_1)&\text{if $q>1$}.
\ecs
\]
\etr

\bcr
The cohomological classes of the following cocycles in $C^*_{(n)}(L_{-1})$
\begin{gather*}
u_{a,b}(n):=e_{2a}\wg e_{2b},\quad\text{where\; $2a+2b=n$\; for integers\; $0\<a<b$},\\
v(n):=\sum_{r=0}^{\frac{n}{4}}e_{\frac{n}{2}-2r-1}\wg e_{\frac{n}{2}+2r+1}\qquad\text{if\; $n\equiv 0\mod 4$}
\end{gather*}
constitute a basis of the space $H^2_{(n)}(L_{-1})$.
In particular,
\[
\dim H^2_{(n)}(L_{-1})=\left\lfloor\frac{n}{4}\right\rfloor+1.
\]
\ecr

Thanks to a well-known interpretation of the $2$-dimensional cohomology
of Lie algebras with trivial coefficients (see \cite{Wei}, Sec.7.6),
the cocycles $u_{a,b}(n)$ and $v(n)$ explicitly describe the basis of the space of equivalence
classes of the central extensions with one-dimensional kernel of the Lie algebra ${W(\mathbb{Z}_2)=L_{-1}}$.



\begin{thebibliography}{12}

\bibitem{MR1634067}
Andrews,~G.~E.,
\newblock{\em The theory of partitions},
\newblock {Cambridge University Press, Cambridge}, 1998.

\bibitem{MR874337}
Fuks,~D.~B.[Fuchs,~D.~B.],
\newblock {\em Cohomology of infinite-dimensional {L}ie algebras}.
\newblock Contemporary Soviet Mathematics. Consultants Bureau, New York, 1986.

\bibitem{MR0440631}
Gel$'$fand,~I.~M.,
\newblock The cohomology of infinite dimensional {L}ie algebras: some questions
  of integral geometry.
\newblock In {\em Actes du Congr\`es International des Math\'ematiciens (Nice,
  1970), Tome 1}, pages 95--111. Gauthier-Villars, Paris, 1971.

\bibitem{MR0339298}
Gon{\v{c}}arova ~L.~V.,
\newblock{Cohomology of {L}ie algebras of formal vector fields on the line},
\newblock \newblock{Functional Analysis and Its Applications}, 1973, 7:2, 6--14; 1973, 7:3, 33--44.

\bibitem{Wei}
Weibel,~C.~A.,
\newblock {\em An Introduction to Homological Algebra},
\newblock {Cambridge studies in advanced mathematics {\bf 38}},
Cambridge University Press 1994.

\bibitem{MR1254731}
Weinstein,~F.~V.,
\newblock {Filtering bases: a tool to compute cohomologies of abstract
subalgebras of the {W}itt algebra},
\newblock In Fuchs,~D.~B.(ed.) {\em Unconventional Lie algebras}, volume~17 of {\em Adv. Soviet
Math.}, pages 155--216. Amer. Math. Soc., Providence, RI, 1993.

\bibitem{W}
Weinstein,~F.~V.,
\newblock{Filtering bases and cohomology of nilpotent
subalgebras of the Witt algebra and algebra of loops in $sl_2$},
\newblock{Functional Analysis and Its Applications}, 2010, 44:1, 4--21

\bibitem{W1}
Weinstein,~F.~V.,
\newblock {Two combinatorial formulas concerning marked partitions},
\newblock \url{arXiv:0805.1467}, 2017.

\end{thebibliography}
\end{document}